\documentclass[12pt]{article}
\usepackage{amsmath}
\usepackage{graphics,graphicx,epsfig}
\usepackage{amssymb}
\usepackage{epstopdf}
\usepackage{sidecap}
\usepackage{appendix}
\DeclareFontFamily{OMX}{yhex}{}
\DeclareFontShape{OMX}{yhex}{m}{n}{<->yhcmex10}{}
\DeclareSymbolFont{yhlargesymbols}{OMX}{yhex}{m}{n}
\DeclareMathAccent{\wideparen}{\mathord}{yhlargesymbols}{"F3}

\newcommand{\Z}{{\mathbb Z}}

\begin{document}

\begin{titlepage}
\vfill
\begin{flushright}
UCPHAS-16-GR05
\end{flushright}

\vfill
\begin{center}
\baselineskip=16pt

{\Large\bf  Golden-Ratio-Based Rectangular Tilings}
\vskip 0.15in
\vskip 10.mm

{\bf  M S Bryant and D W Hobill}

\vskip 0.4cm
Department of Physics and Astronomy, University of Calgary, \\ Calgary, AB, T2N 2W1, Canada \\

\vskip 0.1 in Email: \texttt{mbryant@ucalgary.ca, hobill@ucalgary.ca}

\vspace{6pt}
\end{center}
\vskip 0.2in
\par
\begin{center}{\bf Abstract}
 \end{center}\begin{quote}
A golden-ratio-based rectangular tiling of the first quadrant of the Euclidean plane is constructed by drawing vertical and horizontal grid lines which are located at all even powers of $\phi$ along one axis, and at all odd powers of $\phi$ on the other axis.  The vertices of the rectangles formed by these lines can be connected by rays starting at the origin having slopes that are odd powers of $\phi$.  A refinement of this tiling results in the familiar one with horizontal and vertical grid lines at every power of $\phi$ along each axis. Geometric proofs of the convergence of several known power series' in $\phi$ are provided.

\vfill
\vskip 2.mm
\end{quote}
\hfill
\end{titlepage}

\section{Introduction}

Golden-ratio-based tilings that fill the first quadrant of the Euclidean plane have proven to be of interest, since they lead to some geometric methods for proving certain relationships obeyed by the golden ratio and the Fibonacci numbers.  A tiling pattern introduced and discussed here, provides an alternative to those that have already appeared in the literature \cite{bic,pos}.

The pattern, albeit in a tilted and more skeletal form, first arose in a time series representation of spatial curvature oscillations occurring near the initial cosmological singularity of a Bianchi-IX vacuum spacetime.  The golden ratio and its inverse arose directly from the Einstein equations, and surprisingly, the time series  evolution of the curvature tensor components produced a pattern of self-similar golden rectangles \cite{bry}.  In the present paper, the pattern has been rotated and extended to cover the first quadrant of the Euclidean plane.  It could be subsequently reflected about the fundamental axes to fill the entire plane, but this paper discusses it as a first-quadrant tiling pattern.

The tiling is formed from a grid of horizontal and vertical lines that intersect one axis at odd powers of $\phi$ and the other at even powers of $\phi$.  It is called the AP$\phi$ tiling or the Alternating-Power-of-$\phi$ tiling. The pattern is called alternating because when counting through integer (negative, zero and positive) powers of $\phi$, grid lines appear on one axis, then the other, in an alternating fashion.  By adding extra lines so that each axis has a grid line at every integer power of $\phi$, the tiling can be subdivided to create a previously discussed tiling \cite{bic} (herein called the EP$\phi$ tiling or the Every-Power-of-$\phi$ tiling).

In each of these tilings (AP$\phi$ and EP$\phi$), vertices can be connected by a concurrent family of rays (emanating from the origin) with slopes equal to integer powers of $\phi$.  In this paper, these tilings (and extensions formed by further subdivisions), are used to prove the convergence of some known power series' (powers of $\phi$) and one formula relating $\pi$ and $\phi$.  Two patterns provide visual illustrations of the breeding pattern of Fibonacci's rabbits \cite{fib}.

\section{Some Formulas in $\phi$}

The golden ratio, $\phi=\frac{1+\sqrt{5}}{2}$, has the unique property that $\phi^2=\phi+1$ (or $\phi^{-1}=\phi-1$).  Multiplying through by any power of $\phi$ gives the recursion relation:
\begin{equation}
\label{PhiSplit}\phi^{n+2}=\phi^{n+1}+\phi^{n}\mbox{,}
\end{equation}
which will prove to be useful later.  Four known power series formulas,
% \emph{Golden rectangles} are rectangles which have a length-to-width ratio of $\phi$ and these
%  sometimes play a role in art and design \cite{liv}.
\begin{eqnarray}
\label{formula1} \phi^{n-1} + \phi^{n-3}+\phi^{n-5}+\cdots+ \phi^{n-2k+1}+\cdots & = & \phi^{n}\mbox{,} \\
\label{formula2} \phi^{n-1} + \phi^{n-2}+\phi^{n-3}+\cdots+\phi^{n-k}+\cdots & = & \phi^{n+1}\mbox{,} \\
\label{formula4} 1\phi^{n-1} + 1\phi^{n-3}+2\phi^{n-5}+\cdots+F_k\phi^{n-2k+1}+\cdots & = & \frac{1}{2} \phi^{n+2}\mbox{,} \\
\label{formula3} 1\phi^{n-1} + 2\phi^{n-2}+3\phi^{n-3}+\cdots+k\phi^{n-k}+\cdots & = & \phi^{n+3}\mbox{,}
\end{eqnarray}
where $k=1,2,3,\ldots$ and $F_k$ is the $k^{\textnormal {th}}$ Fibonacci number, will be proven geometrically.  They are generalized versions of those discussed and geometrically proven in \cite{bic}.  The only difference here is that new tilings provide some alternate proofs.

\section{Construction of the Tiling Pattern}

\begin{figure}
\begin{center}
\includegraphics[width=10.0cm]{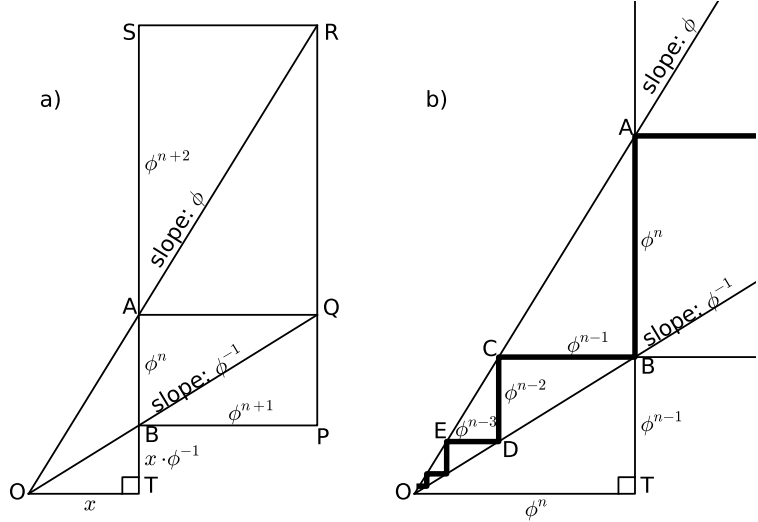}
\caption{Showing construction of a pattern which can be used to create the AP$\phi$ tiling pattern.  The divider is shown as the thick zig-zag path.} \label{Construct}
\end{center}
\end{figure}

Consider the first quadrant of 2D Euclidean space with an as yet unspecified origin.  Construct a landscape-oriented golden rectangle  $ABPQ$, such that point B is closest to the origin and $\overline{BA}$ and $\overline{BP}$ are parallel to the $y$- and $x$-axes respectively and have lengths $\phi^{n}$ and $\phi^{n+1}$ respectively, as shown in Figure \ref{Construct}a.  Extend $\overline{BA}$ and $\overline{PQ}$  in the positive y-direction for a distance  $\phi^{n+2}$ to points $S$ and $R$ respectively, thus forming another golden rectangle $ASRQ$, sharing side $\overline{AQ}$ with $ABPQ$ as shown.

Extend the two rectangle diagonals, $\overline{RA}$ and $\overline{QB}$ until they meet at $O$, as shown.  The diagonals have slopes $\phi$ and $\phi^{-1}$ respectively.  Now designate $O$ as the origin.  Extend $\overline{AB}$ down to point $T$ on the $x$-axis.  $\overline{OT}$ has some fixed but unknown length, $x$, and $\overline{TB}$ has length $x\cdot\phi^{-1}$.  To find $x$, use the slope of $\triangle AOT$'s hypotenuse:
\begin{equation}
\frac{\phi^n+x\phi^{-1}}{x}=\phi\mbox{.} \\
\end{equation}
Solving for $x$, one obtains
\begin{equation}
x=\frac{\phi^n}{\phi-\phi^{-1}}=\frac{\phi^n}{\phi-(\phi-1)}=\phi^{n}\mbox{.} \\
\end{equation}

In Figure \ref{Construct}b, $x$ has been replaced with its now known value.  Also shown in Figure \ref{Construct}b is an infinite zig-zag of horizontal and vertical line segments (hereafter called \emph{the divider}) descending towards the origin between the rays $\overrightarrow{OA}$ and $\overrightarrow{OB}$.   To form the divider, start at $B$ and proceed horizontally leftward to meet $\overrightarrow{OA}$ at $C$, then downward to meet $\overrightarrow{OB}$ at $D$, then leftward to meet $\overrightarrow{OA}$ at $E$, and so on, down to convergence at the origin.  The divider also alternates between $\overrightarrow{OA}$  and $\overrightarrow{OB}$ extending infinitely in the increasing direction. 

Using slope conditions on the right triangles ($\triangle ABC$, $\triangle BCD$, and so on) at the divider,  it is clear that each segment of the divider is a factor of $\phi$ smaller than the one before when proceeding towards the origin. $\overline{OT}$ is the sum of the horizontal divider segments above it, proving (\ref{formula1}).  Adding all the divider segments (horizontal and vertical) from $B$ down to $O$ gives 
\begin{equation}
\phi^{n-1}+\phi^{n-2}+\phi^{n-3}+\cdots=\phi^{n-1}+\phi^{n}=\phi^{n+1}\mbox{,}
\end{equation}
proving (\ref{formula2}). The final step used (\ref{PhiSplit}).

\begin{figure}
\begin{center}
\includegraphics[width=11.0cm]{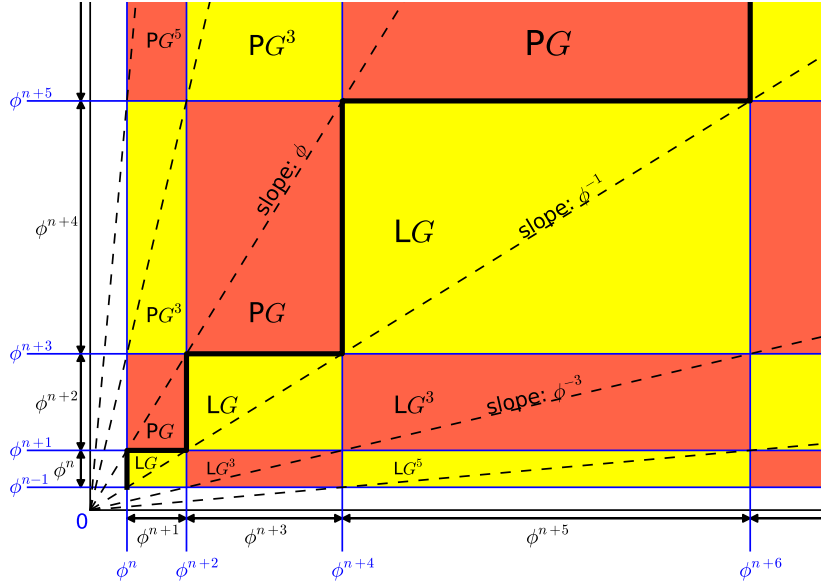}
\caption{The AP$\phi$ tiling pattern (extended from Figure \ref{Construct}) showing power-of-$\phi$ grid lines alternating between the axes and the rays with odd power-of-$\phi$ slopes.  The grid lines converge to the axes where details cannot be shown.} \label{APphi}
\end{center}
\end{figure}

\section{The AP$\phi$ Tiling Pattern}

Extending all divider line segments in Figure \ref{Construct}b to full 1st-quadrant size creates Figure \ref{APphi}, the AP$\phi$ tiling pattern.   In each axis direction, grid lines are spaced every second power apart, as are lines $\overleftrightarrow{BP}$, $\overleftrightarrow{AQ}$ and $\overleftrightarrow{SR}$ in Figure \ref{Construct}a.  This is shown by tick marks and spacings along each axis in Figure \ref{APphi}.  The divider is again shown as a thick zig-zag of line segments.  

It is clear that in the AP$\phi$ tiling, all fundamental rectangles (formed by adjacent grid lines in each dimension), have integer-power-of-$\phi$ dimensions (i.e. are $\phi^j$ by $\phi^i$ rectangles, $i,j \in \Z$).  Nestled underneath the divider is a series of golden rectangles with landscape orientation (labelled as $LG$ for Landscape Golden).  These rectangles can be visualized as rectangular beads, strung by their diagonals, along the ray with slope $\phi^{-1}$, each bead being a factor of $\phi^2$ in linear measure larger than the one before it for progression away from the origin.  Nestled above the divider is a similar set of golden rectangles strung in a similar manner along the ray with slope $\phi$.  In this case they have portrait orientation (and so are labelled PG).  Both sets together represent all  golden rectangles whose side lengths are integer powers of $\phi$.

Outside of these golden rectangles, the rays with slopes $\phi^{-3}$ and $\phi^{3}$ each exhibit the same string-of-beads pattern and together contain all rectangles with power-of-$\phi$ dimensions where the length-to-width ratio is $\phi^{3}$ (labelled appropriately as $LG^3$ and $PG^3$).  This continues with the subsequent rays, which have rectangles with length-to-width ratios of $\phi^{5}$,  $\phi^{7}$, and so on.  Note that all rectangles below the divider have landscape orientation, while those above have portrait orientation.

In the AP$\phi$ tiling, all rectangles with integer-power-of-$\phi$ dimensions for which the length-to-width ratio is  a positive odd power of $\phi$ are represented exactly once, and those are the only ones represented.  The AP$\phi$ pattern is invariant under even-power-of-$\phi$ dilatations about the origin and also under transformations that are a single composition of an odd-power-of-$\phi$ dilatation about the origin with a reflection about the ray $y=x$.

\begin{figure}
\begin{center}
\includegraphics[width=9.0cm]{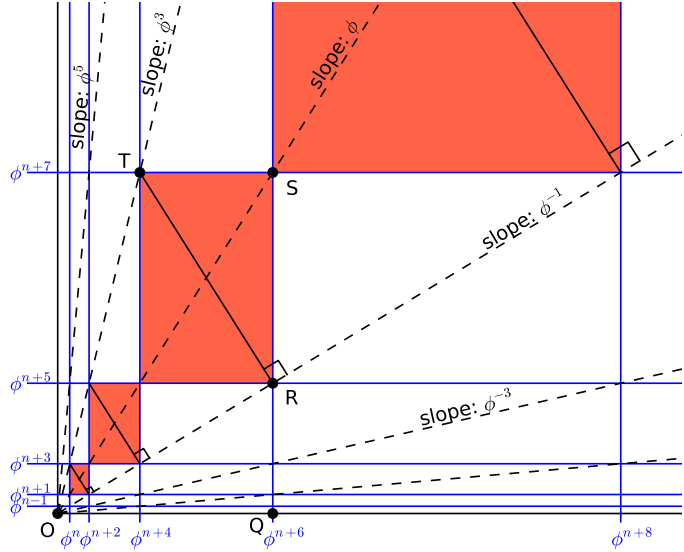}
\caption{The beaded ray of portrait golden rectangles from the AP$\phi$ tiling, with markings and points for proof of a formula relating $\pi$ and $\phi$.} \label{PiDiv4}
\end{center}
\end{figure}

\section{A Formula linking $\pi$ and $\phi$}

In the  AP$\phi$ tiling, the rays with slopes $\phi^3$ and $\phi^{-1}$ are separated by an angle of $45^\circ$, as the following shows.  The AP$\phi$ tiling in Figure \ref{PiDiv4} features the slope-$\phi$ ray whose beads are portrait-oriented golden rectangles.  Consider the rectangle displaying diagonal $\overline{TR}$, which is an arbitrary one since the scaling factor, $n$, is unspecified. $\angle OQR=\angle RST=90^\circ$,  $|\overline{OQ}|=|\overline{RS}|=\phi^{n+6}$ and $|\overline{QR}|=|\overline{ST}|=\phi^{n+5}\Rightarrow \triangle OQR \equiv \triangle RST$ (by the Side-Angle-Side congruence theorem)  $\Rightarrow |\overline{OR}|=|\overline{RT}|$.  $\overline{RT}$ has slope $-\phi$ and $\overline{OR}$ has slope $\phi^{-1}$ so the product of the slopes is $(-\phi) (\phi^{-1})=-1$, so $\overline{OR}\perp\overline{RT}$.  Clearly $\triangle ORT$ is a 45-45-90 triangle and $\angle ROT = 45^\circ = \pi/4$ radians. 

This means $\angle ROT = \angle QOT - \angle QOR$ implies
\begin{equation}
\frac{\pi}{4}=\arctan{(\phi^3)}-\arctan{(\phi^{-1})}
\end{equation}
which is a formula relating $\pi$ and $\phi$.  It is expected that there are geometric proofs for many other formulas hidden in the AP$\phi$ tiling.  The cosmological time series pattern \cite{bry} that was the impetus to investigate the AP$\phi$ tiling was a tilted pattern where the slope-$\phi^{-1}$ ray coincided with the $x$-axis and the slope-$\phi^3$ ray coincided with the ray $y=x$.  The self-similar golden rectangles produced by the time series are the shaded ones in Figure \ref{PiDiv4}.

\section{Division of Power-of-$\phi$-Dimension Rectangles}\label{SectionRDivide}

\begin{figure}
\begin{center}
\includegraphics[width=9.0cm]{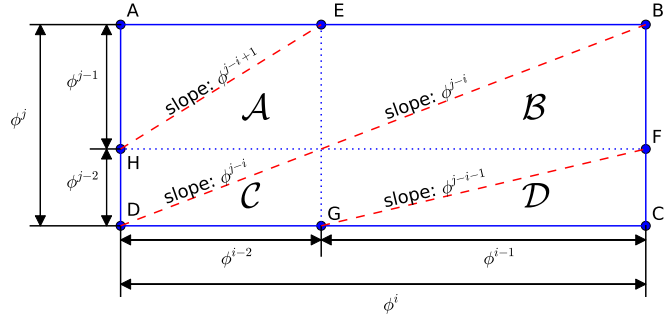}
\caption{An arbitrary integer-power-of-$\phi$ rectangle (measuring $\phi^i$ by $\phi^j$) is divided according to the golden ratio in each dimension, with the smaller portion coming first in each instance.  The sub-rectangles $\mathcal{A}$, $\mathcal{B}$, $\mathcal{C}$ and $\mathcal{D}$, have diagonals with slopes of $\phi^{j-i+1}$, $\phi^{j-i}$, $\phi^{j-i}$ and $\phi^{j-i-1}$, respectively.} \label{RectSubDiv}
\end{center}
\end{figure}

Figure \ref{RectSubDiv} shows a rectangle with arbitrary integer-power-of-$\phi$ dimensions ($\phi^i$ by $\phi^j$) being subdivided using the golden ratio relation, as given in (\ref{PhiSplit}).  If division is performed as shown (the smaller segment preceding the larger one), the rectangle is divided into the four sub-rectangles labelled $\mathcal{A}$, $\mathcal{B}$, $\mathcal{C}$ and $\mathcal{D}$. It is easily seen that the diagonals of $\mathcal{B}$ and $\mathcal{C}$ have the same slope as the parent rectangle, whereas the diagonal of $\mathcal{A}$ ($\mathcal{D}$) has a slope that is $\phi$ times greater (less) than the parent rectangle.

Dividing up rectangles which have integer-power-of-$\phi$ dimensions produces sub-rectangles that also have integer-power-of-$\phi$ dimensions.  If the rectangles are part of a tiling pattern, sub-patterns can be made by repeatedly subdividing rectangles.  In general, the grid locations for the vertices of such sub-rectangles become more complicated due to the creation of sums of various integer powers of $\phi$ that do not conveniently simplify.  Complicated but useful tiling patterns can be made in this way to aid in geometric proofs of the convergence of power series' in $\phi$, as will be shown later in proofs of (\ref{formula4}).

\section{The EP$\phi$ Tiling Pattern}

\begin{figure}
\begin{center}
\includegraphics[width=11.0cm]{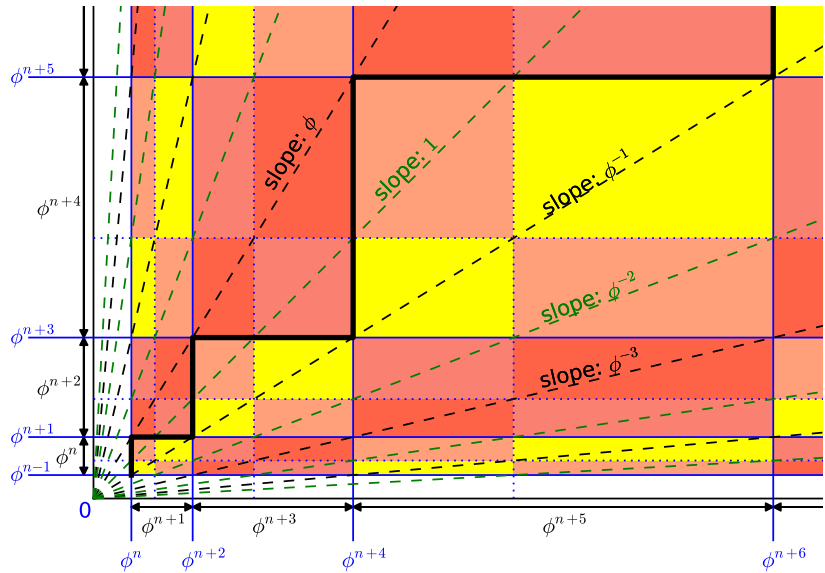}
\caption{The EP$\phi$ tiling pattern has grid lines at every integer power of $\phi$ along each axis.  This is produced by subdivision of the AP$\phi$ tiling (Figure \ref{APphi}) in the manner that Figure \ref{RectSubDiv} illustrates.} \label{EPphi}
\end{center}
\end{figure}

Applying the division shown in Figure \ref{RectSubDiv} to every fundamental rectangle in the AP$\phi$ tiling, produces the EP$\phi$ tiling (see Figure \ref{EPphi}).  It has grid lines located at each power of $\phi$ on each axis and is a pattern which commonly appears in the literature \cite{bic, pos}.

If we consider any two adjacent  rays of the AP$\phi$ tiling (Figure \ref{APphi}), they will have slopes $\phi^k$ and $\phi^{k+2}$ for the lower and upper rays respectively, for some $k$, where $k$ is odd.  Each ray can be visualized as supporting a sequence of self-similar rectangular beads with diagonals positioned along the ray.  These beaded rays fit snugly together, one to the next.

When these Figure \ref{APphi} rectangle/beads are subdivided (in Figure \ref{RectSubDiv} fashion), the $\mathcal{D}$ sub-rectangles from the upper ray beads will have diagonal slopes of $\phi^{k+1}$ (one power less than the slope of the upper ray) and the $\mathcal{A}$ sub-rectangles from the lower ray beads will also have diagonal slopes of $\phi^{k+1}$ (one power more than the slope of the lower ray).  These sub-rectangles combine to form a beaded ray with slope $\phi^{k+1}$.

\begin{figure}
\begin{center}
\includegraphics[width=10.0cm]{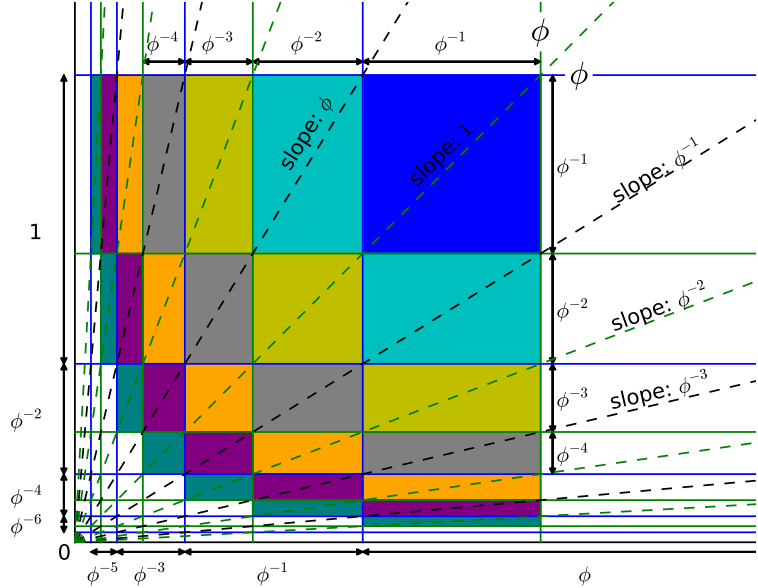}
\caption{Showing the $\phi$-by-$\phi$ square next to the origin of the EP$\phi$ tiling, with shading appropriate for proving an instance of (\ref{formula3}).} \label{SequentCoefficients}
\end{center}
\end{figure}

The EP$\phi$ tiling has the same beads-on-the-ray layout as the AP$\phi$ tiling, except the EP$\phi$ tiling has beaded rays for every power of $\phi$, instead of only for the odd powers.  Further subdivision of the rectangles of the EP$\phi$ tiling leads to $\mathcal{A}$ and $\mathcal{D}$ sub-rectangles with differing slopes (not aligning along a ray), so further subdivision patterns do not have the same beads-on-a-ray structure.

The  EP$\phi$ tiling has reflection symmetry about the ray $y=x$ and is invariant under any integer-power-of-$\phi$ dilatation about the origin.  The beads along the $y=x$ ray consist of all possible squares that have integer-power-of-$\phi$ side-length, each size appearing exactly once.  In the  EP$\phi$ tiling, except for the squares, every other possible rectangle that has integer-powers-of-$\phi$ length and width occurs exactly twice, once in landscape orientation, and once in portrait orientation.

Figure \ref{SequentCoefficients} (also seen in \cite{bic}) shows the EP$\phi$ tiling (with specific powers of $\phi$ on the axes), restricted to the $\phi$-by-$\phi$ section near the origin.  The tiles have been coloured to illustrate a geometric proof of (\ref{formula3}) in the instance $n=-1$:
\begin{equation}
\label{SequentFormula} 1\phi^{-2} + 2\phi^{-3}+3\phi^{-4}+\cdots = \phi^{2}\mbox{.}
\end{equation}
It is sufficient to prove the formula for just one value of $n$, since multiplying through by a power of $\phi$ can then establish it for other $n$ values.  The proof should be clear from the diagram, or one can refer to \cite{pos}.

\section{The Fibonacci Relations}

\begin{figure}
\begin{center}
\includegraphics[width=12.0cm]{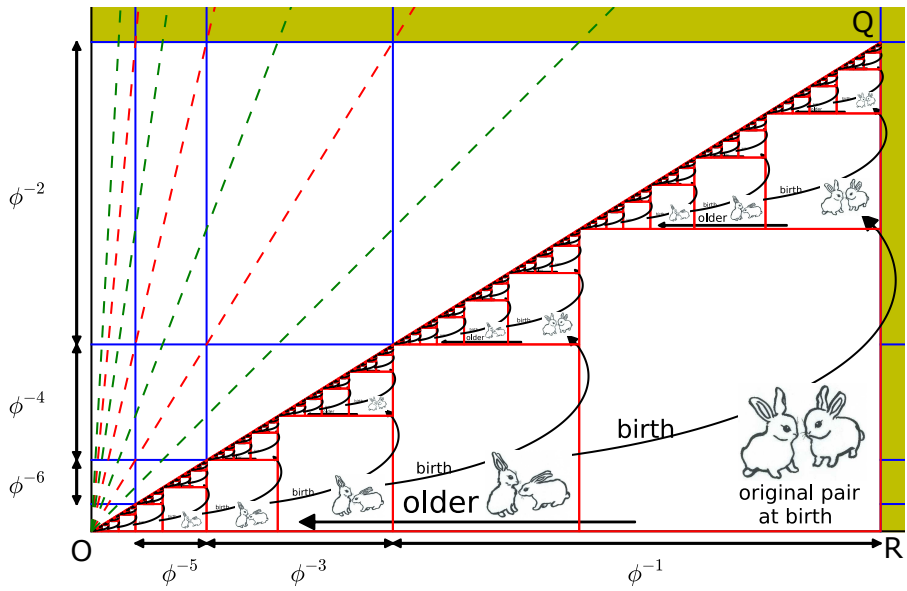}
\caption{A tiling of squares that fills a triangular half of a golden rectangle, geometrically proving an instance of (\ref{formula4}) and providing a visualization of Fibonacci's rabbits.} \label{Rabbtri}
\end{center}
\end{figure}

\begin{figure}
\begin{center}
\includegraphics[width=12.0cm]{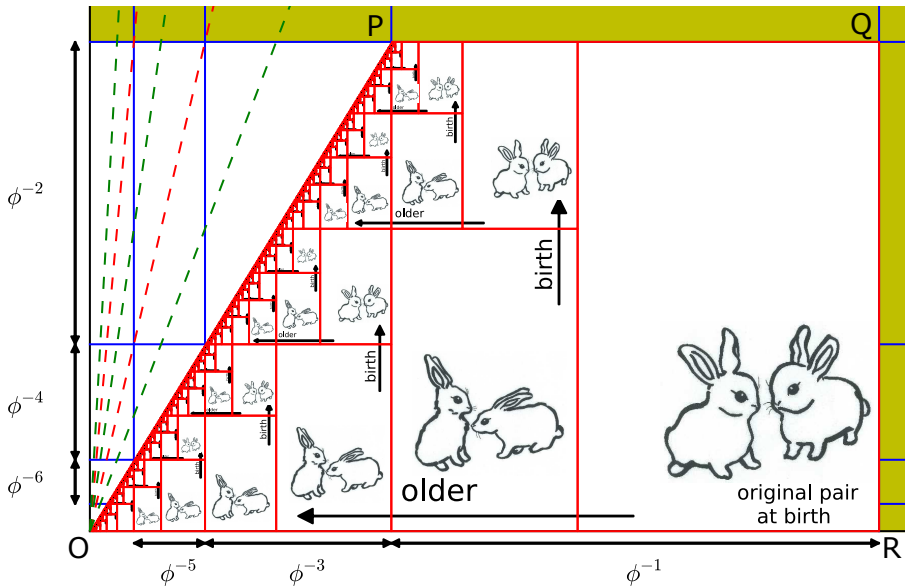}
\caption{A tiling of golden rectangles that fills a trapezoid, geometrically proving an instance of (\ref{formula4}) and providing a visualization of Fibonacci's rabbits.} \label{Rabbzoid}
\end{center}
\end{figure}

As discussed in Section \ref{SectionRDivide}, each tile in an integer-power-of-$\phi$ rectangular tiling pattern can easily be subdivided as many times as desired to form a customized tiling pattern in which all tiles have dimensions which are powers of the golden ratio.  This is done in Figures \ref{Rabbtri} and \ref{Rabbzoid}, where self-similar patterns have been created to prove instances of (\ref{formula4}).  For each of these figures, different tile sizes represent different terms of the series.

To show that the number of tiles for each tile size is the correct Fibonacci number, Fibonacci's original concept of counting rabbit pairs has been employed \cite{fib}.  In that paper, the Fibonacci numbers are the number of rabbit pairs in existence each month, where the first month has only a single pair at birth.  Newly born pairs do not reproduce for their first month, but produce one more pair (capable of future breeding) in each month after that.  It is assumed that all rabbits live forever.

In each Figure, the largest tile (labelled ``original pair at birth"), represents the rabbit pair that begins the breeding pattern at the start of Month 1.  The rabbit pairs in each subsequent month are represented by tiles which are a factor of $\phi$ smaller than those of the month before.  Since there is a one-to-one correspondence between rabbit pairs in the month and tiles of that month's size, in Month $n$, the total number of tiles of the appropriate size will be $F_n$, representing the number of rabbit pairs present in that month (counting all pre-existing and newly-born).  The fact that each rabbit pair lives forever is represented by the series of tiles (lower edges aligned) to the left of each baby pair.

Figure \ref{Rabbtri} illustrates the $n=-3$ instance of (\ref{formula4}):
\begin{equation}
\label{RabbtriFormula} 1\phi^{-4} + 1\phi^{-6}+2\phi^{-8}+\cdots+F_k\phi^{-2k-2}+\cdots = \frac{1}{2} \phi^{-1}\mbox{.}
\end{equation}
The tiles are all squares and correspond to rabbit pairs.  For Month $k$ each square has area $\phi^{-2k-2}$ and there are $F_k$ of them so the total area of Month $k$ squares is  $F_k\phi^{-2k-2}$ which is the $k$th term of  the LHS of (\ref{RabbtriFormula}) so the LHS of (\ref{RabbtriFormula}), is the total area of all rabbit pair tiles in the Figure \ref{Rabbtri} tiling.

Each of the square tiles in the pattern touches the $\triangle OQR$ hypotenuse, $\overline{OQ}$.  As the breeding pattern progresses, the new points where $\overline{OQ}$ is touched are always proportionally spaced, so it is clear that the tiling pattern will, in its limit, completely fill $\triangle OQR$.  $\triangle OQR$ has area $\frac{1}{2}|\overline{OR}||\overline{RQ}| = \frac{1}{2}(1)(\phi^{-1})=\frac{1}{2}\phi^{-1}$ which is the RHS of (\ref{RabbtriFormula}), completing the proof.  The Figure \ref{Rabbtri} pattern appears in \cite{hut}, without mathematical analysis.

Figure \ref{Rabbzoid} illustrates the $n=-2$ instance of (\ref{formula4}):
\begin{equation}
\label{RabbzoidFormula} 1\phi^{-3} + 1\phi^{-5}+2\phi^{-7}+\cdots+F_k\phi^{-2k-1}+\cdots = \frac{1}{2}\mbox{.}
\end{equation}
The tiles are all golden rectangles and correspond to rabbit pairs.  For Month $k$ each rectangle has area $\phi^{-2k-1}$ and there are $F_k$ of them so the total area of Month $k$ rectangles is  $F_k\phi^{-2k-1}$ which is the $k$th term of  the LHS of (\ref{RabbzoidFormula}), so the LHS of (\ref{RabbzoidFormula}) is the total area of all rabbit pair tiles in the Figure \ref{Rabbzoid} tiling.

Trapezoid $OPQR$ has area $\frac{1}{2}(|\overline{PQ}|+|\overline{OR}|)|\overline{RQ}| = \frac{1}{2}(1+\phi^{-1})(\phi^{-1})=\frac{1}{2}(\phi)(\phi^{-1}) = \frac{1}{2}$ which is the RHS of (\ref{RabbzoidFormula}).  It remains to prove that the tiling pattern completely fills trapezoid $OPQR$.

There is a sequence of baby-rabbit-pair tiles, where each tile in the sequence touches $\overline{PQ}$, the top side of the trapezoid.  The sum of the widths of these tiles is $\phi^{-2}+\phi^{-4}+\phi^{-6}+\cdots=\phi^{-1}=|\overline{PQ}|$, ((\ref{formula1}) with $n=0$ was used), so the pattern fills the trapezoid along the top side, $\overline{PQ}$.

The sum of the widths of the tiles corresponding to the \emph{original rabbit pair} existing for every month is $\phi^{-2}+\phi^{-3}+\phi^{-4}+\cdots=\phi^{0}=1=|\overline{OR}|$, ((\ref{formula2}) with $n=-1$ was used), so the series of tiles representing the \emph{original pair} at all ages converges at $O$, meaning the pattern fills the base of the trapezoid, $\overline{OR}$.

It remains to show that the tiling pattern fills the trapezoid up to the boundary, $\overline{OP}$.  For each \emph{other rabbit pair} (i.e. not the original pair), there is also a series of tiles representing them at all ages.  It starts at the lower right corner of the tile representing them as babies and extends leftward.  The location of that initial corner depends on events that happened before the pair were born.

The convergence point for a rabbit pair series is the same as the convergence point, $O$, of the original pair series, except that every time one of the rabbit pairs' ancestors gives birth, there is (1) a rightward shift by\begin{em} the \begin{bf}width\end{bf} of the ancestor tile for the month when the birth happened\end{em}, and (2) an upward shift by\begin{em} the \begin{bf}height\end{bf} of the ancestor tile for the month when the birth happened\end{em}.  The ratio of \emph{upward} shift to \emph{rightward} shift is always the golden ratio.

The series convergence point for \emph{any rabbit pair} is shifted \emph{upward} and \emph{rightward} along the ray $y=\phi x$ whenever one of the rabbit pairs' ancestors gives birth.  But this ray is $\overline{OP}$.  Therefore, the series for every rabbit pair always converges to a point on $\overline{OP}$, at the same height as the lower right corner of the tile representing the pair at birth.  The vertical locations of the convergence points on $\overline{OP}$ are always proportionally spaced, so it is clear that the tiling pattern will, in its limit, completely fill trapezoid $OPQR$, so the proof is complete. 

\section{Conclusion}

Golden-ratio-based tiling patterns have been of interest from both an aesthetic and a mathematical point of view.  The main pattern that has appeared in the literature (referred to as the EP$\phi$ tiling) is created by orthogonal grid lines appearing at every power of $\phi$ along both axes of the first quadrant of a Cartesian coordinate system.  The AP$\phi$ tiling, introduced here, has grid lines at all powers of phi, but they alternate between the two orthogonal axes, so one axis contains all the even power-of-$\phi$ locations and the other the odd.  This results in a tiling where each tile is unique.  The EP$\phi$ tiling is obtained from the AP$\phi$ tiling by sub-division of the AP$\phi$ tiles using the fundamental relation $\phi^2 = \phi +1$ (which is merely adding extra grid lines).  The EP$\phi$ tiling has pairs of congruent tiles (one of each pair having landscape orientation and the other portrait), except for the set of self-similar $\phi^i \times \phi^i$ squares where each square appears only once.

This paper has noted that each of the two tilings has a discrete family of rays emanating from the origin and passing through the diagonals of the rectangular tiles.  Each ray can be associated with a self-similar ``rectangular-beads-on-a-ray'' structure.  The EP$\phi$ tiling consists of beaded-rays that have slopes of every integer power of $\phi$, while the AP$\phi$ tiling consists of beaded-rays with slopes of every odd integer power of $\phi$.

The EP$\phi$ tiling pattern is useful for providing geometric proofs of the convergence of certain power series' of $\phi$.  The AP$\phi$ tiling's sparser structure (missing grid lines compared to the EP$\phi$ tiling) serves to highlight other properties and proofs.  It facilitated recognition of the ``cosmological curvature pattern'', leading to a relation between $\pi$ and $\phi$. By stepping back from the EP$\phi$ tiling and looking more generally at tile subdivisions, two self-similar tilings were obtained which aided visualization of the pattern of Fibonacci's breeding rabbit model and provided proofs of the convergence of a power series of $\phi$.

It can be expected that other relationships remain hidden in the tilings, dividers and families of beaded-rays introduced here.  It is hoped that the perspectives discussed in this manuscript will inspire others in the search for alternative mathematical relations that might have a basis in geometrically self-similar golden ratio patterns.

\section{Acknowledgments}

The authors wish to thank Holly Farrell (step-daughter of MB), for creating the rabbit pair images.  One of us (DH) gratefully acknowledges financial support from the National Sciences and Engineering Research Council of Canada (NSERC).

\end{document}